\RequirePackage{fix-cm}
\documentclass[smallextended]{svjour3}       
\smartqed  
\usepackage{amssymb,amsmath,multirow}
\usepackage{xcolor}

\newtheorem{thm}{Theorem}
\newtheorem{ass}{Assumption}
\newcommand{\vectwo}[2]{\left(\begin{array}{c}#1 \\ #2 \end{array}\right)}
\def\pn{\par\smallskip\noindent}
\def\it proof{\pn {Proof.} }

\begin{document}
\title{A survey of hidden convex optimization
\thanks{This research was supported  by National Natural Science Foundation of China
under grants 11822103, 11571029 and Beijing Natural Science Foundation Z180005.}
}

\titlerunning{Hidden convex optimization}        

\author{Yong Xia}


\institute{
Y. Xia  \at
              LMIB of the Ministry of Education; 
              School of
Mathematics and System Sciences, Beihang University, Beijing,
100191, P. R. China
              \email{ yxia@buaa.edu.cn}
 }

\date{Received: date / Accepted: date}

\maketitle

\begin{abstract}
Motivated by the fact that
not all nonconvex optimization problems are difficult to solve, we survey in this paper three widely-used ways to reveal the hidden convex structure for different classes of nonconvex optimization problems. Finally, ten  open problems are raised.
\end{abstract}

\section{Introduction}

Nowadays, convex programming becomes very popular due to not only  its various real applications but also the charming property  that any local solution also remains global optimality.

However, it does not mean that convex programming problems are easy to solve.   First, it may be difficult to identify convexity. Actually,
deciding whether a quartic  polynomials  is globally convex is NP-hard
\cite{Ah13}.
Second, evaluating a convex function is also not always easy.  For example,
the induced matrix norm
\[
\|A\|_p=\sup_{x\neq 0}\frac{\|Ax\|_p}{\|x\|_p}
\]
is a convex function. But evaluating  $\|A\|_p$ is NP-hard  if $p\neq1,2,\infty$ \cite{H09}.
Third, convex programming problems may be difficult to solve.
It has been shown in \cite{Bu09} that  the general mixed-binary quadratic optimization problems can be equivalently reformulated as convex programming over the copositive cone
\[
{\rm Co}:=\{A\in\Bbb R^{n\times n}:~A=A^T,~x^TAx\ge 0,~\forall  x\ge0, x\in\Bbb R^n \}.
\]
Notice that checking whether $A\in {\rm Co}$ is NP-Complete \cite{M87}.

Like a coin has two sides, there are quite a few nonconvex optimization problems could be globally and efficiently solved in polynomial time. The reason behind this observation is that most of them belong  to the \textit{hidden convex} optimization, i.e., they admit equivalent convex programming reformulations. In the past two decades, hidden convex optimization problems were studied case by case in literature.
In this survey, we summarize three class of approaches to reveal the  hidden convex structure, namely, nonlinear transformation, Lagrangian dual and the primal tight convex relaxation.

The remainder of this paper is organized as follows.
Section 2 presents the most common used nonlinear transformations for different classes of problems.
Section 3 shows how  Lagrangian duality and its variations achieve zero gap. Section 4 summarizes some classes of nonconvex optimization problems admitting tight primal relaxations. As a concluding remark, ten open problems are raised in Section 5.

Throughout the paper, let $v(\cdot)$ be the optimal value of problem $(\cdot)$.
$\Bbb R^n_+$ denotes the nonnegative orthant in $\Bbb R^n$.
For a matrix $A$,
denote by $A\succ(\succeq)0$ that $A$ is positive (semi)definite. The trace of  $A$ is defined as the sum of its diagonal elements, i.e., ${\rm tr}(A)=\sum_{i=1}^nA_{ii}$.
$\lambda_{\min}(A)$ and $\lambda_{\max}(A)$ denote the minimal and maximal eigenvalues of $A$, respectively. $I_n$ is the identity matrix of order $n$.
 $\|x\|=\sqrt{x^Tx}$ stands for the $\ell_2$ norm of a vector $x$. ${\rm diag(x_1,\ldots,x_n)}$ returns a diagonal matrix with diagonal elements $x_1,\ldots,x_n$. $e=(1,1,\ldots,1)^T\in\Bbb R^n$.
Denote by ${\rm conv}\{\Omega\}$ the convex hull of $\Omega$. $|E|$ denotes the number of elements in the set $E$. For a real value $x$,  [$x$] returns the largest integer less than or equal to $x$.

\section{Nonlinear transformation}
In this section, we survey some nonlinear transformation approaches widely used in convexifying different classes of nonconvex optimization problems.

\subsection{A univariate example}
Univariate examples are not  always as simple as  converting the concave function $\sqrt{x}$  to $y(\ge0)$ by
introducing the one-to-one mapping $x=y^2$ with $y\ge 0$.

The approach \cite{X11} for reducing the  duality gap for box constrained nonconvex quadratic
program requires  solving the following nonconvex subproblem
\begin{equation}
({\rm G})~~\max_{\theta} \phi(\theta):=\min\{a_1\theta,a_2\delta^2(\theta)\}, \label{G}
\end{equation}
where $a_1$, $a_2$ are two positive scalars,
\begin{eqnarray}
\delta^2(\theta)&=& \min\sum_{i\in I\cup J}(x_i-y_i)^2\nonumber\\
 &&{\rm s.t.}~x\in C,~ y_i\in[-1, 1],~i\in I,\nonumber\\
 &&~~~~~\sqrt{1-\theta}\le \omega_iy_i\le1,~i\in J,\nonumber
\end{eqnarray}
$C$ is a linear manifold and $\omega_i\in\{-1,1\}$ is given. The objective function $\phi(\theta)$ (\ref{G}) is not concave and thus (G) is a nonconvex minimization.
 Introducing $\sigma=\sqrt{1-\theta}$, one can reformulate $({\rm G})$ as the following convex program:
\begin{eqnarray}
&\max& \sigma   \nonumber\\
&{\rm s.t.}&\left\|\vectwo{\sqrt{a_2}\cdot (x-y)}{\sqrt{a_1}\cdot
\sigma}\right\|\le \sqrt{a_1}, \nonumber\\
&&-1\leq y_i\leq 1,~~~~i\in I,\nonumber \\
&& \sigma\leq \omega_i y_i\leq 1,~~~i\in J,\nonumber \\
&&x\in C, ~ \sigma\geq 0.\nonumber
\end{eqnarray}

\subsection{$p$-th power approach}
Consider the following constrained nonconvex  optimization problems:
\begin{eqnarray*}
{\rm(P)}~~ \min\left\{ f(x)~: g_i(x)\le b_i,i=1,\ldots,m,~x\in S\right\},
\end{eqnarray*}
where $f(x)$, $g_i(x)$ ($i=1,2,\ldots,m$) are strictly  positive over $S$ and $S$ is closed, bounded, connected and of a full dimension.  It is known \cite{L70} that if the  perturbation  function
\[
w(y):=\min_{x\in S}\left\{ f(x)~: g_i(x)\le y_i,i=1,\ldots,m,~x\in S\right\}
\]
is locally
convex in a closed and non-degenerate neighbourhood of $b$, then there is no duality gap between $(P)$
and its Lagrangian dual in the sense that the inner minimization is local
and the outer maximization is in the neighbourhood of the optimal
Lagrangian multiplier.
In order to achieve such a zero duality gap, Li
\cite{Li95} first introduced the  $p$-th power transformation:
\begin{eqnarray*}
({\rm{P}}_p)~~ \min \left\{ f^p(x)~:
  g^p_i(x)\le b^p_i,i=1,\ldots,m,~x\in S\right\}.
\end{eqnarray*}
Under some additional assumptions, Li \cite{Li95}
showed that, for a sufficiently large $p$, the perturbation function of $({\rm{P}}_p)$ is a
convex function of $y^p$ for any $y$ in a neighbourhood of $b$.

For more applications of
the $p$-th power formulation, we refer to \cite{L06,L07} and references therein.
More generalizations of the $p$-th power convexification approach are further studied in  \cite{L05,Wu07}.

Recently, a novel shifted $p$-th power reformulation is introduced in \cite{X16}:
\begin{eqnarray*}
({\rm{P}}_{p,\mu})~~ \min\left\{ f(x)~:
  (g_i(x)+\mu_i)^p\le (b_i+\mu_i)^p,i=1,\ldots,m,~x\in S\right\}.
\end{eqnarray*}
Surprisingly, under the same assumptions, there is a parametric vector
$\mu\in \mathbb{R}^m$ such that $p=3$  is sufficient to guarantee the convexity of the perturbation function of $({\rm{P}}_{p,\mu})$ in terms of $(y+\mu)^p$ with $y$ lying in a neighbourhood of $b$. Besides, the assumption on the strict positiveness of $f(x)$ and $g_i(x)$ is redundant in this new reformulation. For more details, we refer to \cite{X16}.

\subsection{Minimal-volumn ellipsoid cover}
Consider the geometric problem of finding $n$-dimensional ellipsoid of minimal volume covering a set of $m$ given points $a_i,~i=1,\ldots,m$:
\begin{eqnarray}
({\rm MVE})&\min_{Q,x}& u_n\sqrt{{\rm det}(Q^{-1})}\nonumber\\
&{\rm s.t.}& (x-a_i)^TQ(x-a_i)\le1,i=1,\ldots,m,\nonumber
\end{eqnarray}
where $u_n$ is the volume of the unit ball in $\Bbb R^n$ and the objective is the $n$-dimensional volume.

By first introducing $M=Q^{\frac{1}{2}}$ (which is well defined as $Q\succ 0$) and $z=Mx$, and then taking
a logarithmic transformation for the objective function,
we can reformulate  the nonconvex optimization problem (MVE)  as the following convex program \cite{S04}:
\begin{eqnarray}
&\min_{M,z}& -u_n\log {\rm det}(M)\nonumber\\
&{\rm s.t.}& (z-Ma_i)^T(z-Ma_i)\le1,i=1,\ldots,m,\nonumber
\end{eqnarray}
where the fact det$(Q^{-1})=($det$(M))^{-2}$ is used.

\subsection{Multiplicative programming}
Consider the linear multiplicative programming \cite{K95}:
\begin{eqnarray}
{\rm (LMP)}~~&\max_{x\in C}& (d_1^Tx+c_1)(d_2^Tx+c_2) \nonumber\\
&{\rm s.t.}&  d_i^Tx+c_i\ge 0,~ i=1,2, \nonumber
\end{eqnarray}
where $C$ is a convex set. (LMP) is NP-hard \cite{M96} to solve if replacing the maximization with minimization. The objective function of (LMP) is nonconcave but quasiconcave \cite{K95}.  The hidden convexity of (LMP) is viewed by the equivalent convex programming problem in the sense that both sharing the same optimal solution:
\begin{eqnarray}
 &\max_{x\in C}& \log(d_1^Tx+c_1)+\log(d_2^Tx+c_2) \nonumber\\
&{\rm s.t.}&  d_i^Tx+c_i\ge 0,~ i=1,2. \nonumber
\end{eqnarray}
Moreover,
we notice  that (LMP) has the following new second-order cone programming representation:
\begin{eqnarray}
&\max_{x\in C}&~ z \nonumber\\
&{\rm s.t.}&  d_i^Tx+c_i\ge 0,~ i=1,2, \nonumber\\
&& s\le d_1^Tx+c_1,~t\le  d_2^Tx+c_2,  \nonumber\\
&&  \left\|\vectwo{s-t}{2z}\right\|\le  s+t.\nonumber
\end{eqnarray}

\subsection{Geometric programming}
Geometric program was first introduced in
the book \cite{D67}. It is an optimization problem of the form
\begin{eqnarray}
{\rm (GP)}~~&\min& f_0(x) \nonumber\\
&{\rm s.t.}& f_i(x)\le 1,~i=1,\ldots,m,\label{g:1}\\
&& g_j(x)=1,~j=1,\ldots,p,\label{g:2}
\end{eqnarray}
where $g_j~(j=1,\ldots,p)$ are monomials, i.e. the form $c_j\prod_{i=1}^nx_i^{a_{ij}}$ with $c_j>0$, and $f_i~(i=0,\ldots,m)$ are posynomials, i.e., the sum of several monomials. We also assume that the constraints (\ref{g:1})-(\ref{g:2}) implicitly imply that all the variables are positive.

Obviously, the general (GP) is a nonconvex optimization problem. One can verify that   introducing logarithmic change of variables $y_i=\log(x_i)$ and logarithmic transformations of $f_i(x)$ yields a    convex programming reformulation of (GP):
\begin{eqnarray}
 &\min& \log f_0(e^y) \\
&{\rm s.t.}& \log f_i(e^y)\le 0,~i=1,\ldots,m,\\
&& \log g_j(e^y)=0,~j=1,\ldots,p,
\end{eqnarray}
where $(e^y)_i=e^{y_i}$.

\subsection{Fractional programming}
Consider the linear fractional program
\begin{eqnarray}
{\rm (LF)}~~&\min& f(x)= \frac{c^Tx+\alpha}{d^Tx+\beta} \label{Lfrac}\\
&{\rm s.t.}& x\in\Omega:=\left\{x:~ Ax\le b,~x\ge 0\right\}.
\end{eqnarray}
To guarantee that (LF) is well defined, we  assume $\min_{x\in\Omega}d^Tx+\beta>0$. It is trivial to verify that the objective function $f(x)$ (\ref{Lfrac}) is quasi-convex and generally nonconvex. Introducing new variables $y$ and $t$ to replace $\frac{x}{d^Tx+\beta}$ and $\frac{1}{d^Tx+\beta}$, respectively, we can equivalently reformulate (LF) as the following linear programming:
\begin{eqnarray*}
&\min&  c^Ty+\alpha t\\
&{\rm s.t.}&  Ay\le tb,~y\ge 0.
\end{eqnarray*}
The above approach was initially suggested by Charnes and Cooper \cite{C62}. It was then  extended by Schaible \cite{S74} to nonlinear convex fractional programs:
\begin{eqnarray}
{\rm (NLF)}~~&\min&  \frac{g_0(x)}{h(x)}\label{nlf:1}\\
&{\rm s.t.}&  g_i(x)\le 0,~i=1,\ldots,m,\label{nlf:2}
\end{eqnarray}
where $g_i(x)$ ($i=0,1,\ldots,m$) are all convex functions, $h(x)$ is a concave function and $h(x)>0$ over the feasible region. Then, introducing
$y=\frac{x}{h(x)}$ and $t=\frac{1}{h(x)}$ reduces the nonconvex optimization problem (NLF) to the following convex programming:
\begin{eqnarray}
{\rm (CLF)}~~&\min&  t\cdot g_0\left(\frac{y}{t}\right)\label{nlf:3} \\
&{\rm s.t.}&  t\cdot g_i\left(\frac{y}{t}\right)\le 0,~i=1,\ldots,m, \label{nlf:4}\\
&& t\cdot h\left(\frac{y}{t}\right) \le 1,~t>0.\label{nlf:5}
\end{eqnarray}

Minimizing the sum of a linear and a linear fractional function over a polyhedral set is generally NP-hard \cite{M96}. Necessary and sufficient condition for the pseudoconcavity of the objective function is established in \cite{C08}. The pseudoconvex but nonconvex case is the following:
\begin{eqnarray*}
{\rm(SLF)}~~&\min& \frac{a_{1}^{T}x+b_{1}}{a_{2}^{T}x+b_{2}}+a_{2}^{T}x \nonumber  \\
&{\rm s.t.}& Ax=b,~x\ge 0.
\end{eqnarray*}
Based on the Charnes-Cooper transformation, it was shown in \cite{FG2016} that (SLF) enjoys a second-order cone programming reformulation:
\begin{eqnarray}
{\rm (SOCP_{LP})}~~ &\min_{s,t,y,z}& t \nonumber  \\
&{\rm s.t.}& Ay-sb=0, ~y \ge 0,\nonumber \\
&& a_{2}^{T}y+b_{2}s=1, \nonumber \\
&& z=b_2a_{2}^{T}y+b_2^2+a_{1}^{T}y+(b_2^2+b_{1})s, \nonumber \\
&& \|(2a_{2}^{T}y, s-t+z)^T\|\le s+t-z. \nonumber
\end{eqnarray}


\subsection{Eigenvalue and trust-region subproblem}

Define the Rayleigh quotient of a symmetric matrix $A\in \Bbb R^{n\times n}$ as:
\[
R_A(x)=\frac{x^TAx}{x^Tx},~x\neq 0.
\]
The Rayleigh-Ritz formula is well known for the maximal eigenvalue of $A$:
\begin{equation}
\sup_{x\neq 0}R_A(x) = \max_{x^Tx=1} x^TAx. \label{eig}
\end{equation}
Let $A=UDU^T$ be the eigenvalue decomposition of $A$, where $U$ is orthogonal and $D={\rm diag}(\lambda_1,\ldots,\lambda_n)$ with $\lambda_1\le \ldots \le\lambda_n$ being the $n$ eigenvalues of $A$. Then, we have
\begin{eqnarray}
\max_{x^Tx=1} x^TAx&=&\max_{y^Ty=1}y^TDy
 = \max_{\sum_{i=1}^n y_i^2=1}\sum_{i=1}^n \lambda_i y_i^2
  =  \max_{e^Tz=1,z\ge0,\forall i}\sum_{i=1}^n \lambda_i z_i, \nonumber
\end{eqnarray}
where $y=U^Tx$ is replaced in the first equality and it holds that $\|y\|=\|x\|$ as $U$ is orthogonal, the nonlinear transformation $z_i=y_i^2$, $i=1,\ldots,n$ are introduced in the last equality and we note that these mappings are no longer one-to-one.  Problem on the right-hand side of the above chain equalities is a linear programming over a standard simplex, which serves as the hidden convex reformulation of  (\ref{eig}).

Now we move to the well-known Kantorovich inequality, which is related to the following nonconvex optimization:
\[
{\rm (KI)}~~\max_{x\neq 0}\frac{(x^TAx)(x^TA^{-1}x)}{(x^Tx)^2},
\]
where $A\succ 0$.  Based on the same transformation as above, (KI) is reduced to
\begin{eqnarray}
\max_{e^Tz=1,z\ge0}\left(\sum_{i=1}^n \lambda_i z_i\right)
\left(\sum_{i=1}^n \lambda_i^{-1} z_i\right),\nonumber
\end{eqnarray}
which is further equivalent (in the sense that both optimal solutions are the same) to solving
\begin{eqnarray}
&&\max_{e^Tz=1,z\ge0}\min_{t\in [t_1,t_2]} ~t\left(\sum_{i=1}^n \lambda_i z_i\right)+
t^{-1}\left(\sum_{i=1}^n \lambda_i^{-1} z_i\right) \nonumber \\
&&=\min_{t\in [t_1,t_2]}\max_{e^Tz=1,z\ge0} ~t\left(\sum_{i=1}^n \lambda_i z_i\right)+
t^{-1}\left(\sum_{i=1}^n \lambda_i^{-1} z_i\right), \label{t12} \\
&&=\min_{t\in [t_1,t_2]}~ \left\{f(t):=\max_{i=1,\ldots,n} t \lambda_i +
t^{-1} \lambda_i^{-1}\right\}, \label{t123}
\end{eqnarray}
where $t_1=\sqrt{\lambda_1/\lambda_n}$, $t_2=t_1^{-1}$, and (\ref{t12}) follows from the classical von Neumann's minimax theorem. Problem (\ref{t123}) is a univariate convex optimization problem and can be explicitly solved.

The inhomogeneous extension of (\ref{eig}) is the well-known trust-region subproblem (TRS) \cite{G81}:
\begin{eqnarray}
{\rm(TRS)}~~\min \left\{ x^TAx+b^Tx:~ x^Tx\le \delta\right\}, \label{trs}
\end{eqnarray}
where $\delta$ is a positive scalar.
(TRS) (\ref{trs}) plays a great role in the trust-region method \cite{Y15} and also has some other applications such as the constrained least squares \cite{G91}.  In the case that $A\not\succeq 0$,  (TRS) is a nonconvex optimization. Interestingly, it was shown in \cite{M94} that (TRS) has at most one local non-global minimizer.

As above, let $A=U{\rm diag}(\lambda_1,\ldots,\lambda_n)U^T$ be the eigenvalue decomposition of $A$ and  $\widetilde{b}=U^Tb$. By introducing $y=U^Tx$, (TRS) is equivalent to
\begin{eqnarray}
 \min \left\{\sum_{i=1}^n \lambda_iy_i^2+\widetilde{b}_iy_i:~ \sum_{i=1}^ny_i^2\le \delta\right\}.\label{T:1}
\end{eqnarray}
For $i=1,\ldots,n$, by further introducing
\[
y_i=\left\{ \begin{array}{ll}\sqrt{z_i},& \widetilde{b}_i\le 0,\\
-\sqrt{z_i},& \widetilde{b}_i>0, \end{array}\right.
\]
we can reduce (TRS) (\ref{T:1}) to the following convex programming over a  simplex:
\begin{eqnarray*}
 \min \left\{ \sum_{i=1}^n \lambda_iz_i -|\widetilde{b}_i|\sqrt{z_i}:~
 \sum_{i=1}^nz_i\le \delta,~z_i\ge 0,~i=1,\ldots,n\right\}.
\end{eqnarray*}
The above convexification approach is further extended by Ben-Tal and Teboulle \cite{B96} to the two-sided trust-region subproblem \cite{S95}. Similar application in convexifying the following regularized problem
\[
\min_{x\in\Bbb R^n} ~x^TAx+b^Tx +\rho \|x\|^p
\]
with $p>2$ can be found in \cite{H17}.

\subsection{Quadratic matrix programming over orthogonal constraints}

For vectors $u$ and $v$, define the minimal product as
\begin{eqnarray}
\langle u,\ v\rangle_{-}\ :=\min_{\pi}\sum_{i=1}^{n}u_{i}v_{\pi(i)},
\end{eqnarray}
where $\pi$ is a permutation of $1,2,\ldots,n$. Notice that
$\langle u,\ v\rangle_{-}$ is easy to solve by first sorting $u$ and $v$, respectively. $\langle
u,\ v\rangle_{+}$ is similarly defined and computed:
\begin{eqnarray}
\langle u,\ v\rangle_{+}\ :=\max_{\pi}\sum_{i=1}^{n}u_{i}v_{\pi(i)}.
\end{eqnarray}

The quadratic assignment problem (QAP) is a classical combinatorial optimization problem \cite{L97}.  The trace formulation reads as follows:
\begin{eqnarray}
{\rm(QAP)}~~&\min& {\rm tr} (AXB^TX^T) \label{qap:tr:1}\\
&& X\in {\rm\Pi}_n:=\{X\in \{0,1\}^{n\times n}:
Xe=X^Te=e\},\label{qap:tr:2}
\end{eqnarray}
where $A, B\in \Bbb R^{n\times n}$ correspond to flow matrix
and distance matrix in a facility location application,
respectively, ${\rm\Pi}_{n}$ is the set of all $n\times n$
permutation matrices.  The orthogonal relaxation of (QAP) was first proposed  in \cite{F87}:
\begin{eqnarray}
{\rm(O)}~~&\min& {\rm tr} (AXB^TX^T)  \label{O10}\\
&& X\in\mathcal{O}_n:=\{X\in \Bbb R^{n\times n}:\ X^TX =I_n\}. \label{O1}
\end{eqnarray}
Let $\lambda$ and $\mu$ be the vectors composed by the eigenvalues of $A$ and $B$, respectively, i.e.,
\[
A=U{\rm diag}(\lambda)U^T,~B=V{\rm diag}(\mu)V^T.
\]
It is not difficult to show the following result.
\begin{thm}[\cite{F87,R92}]
For any $X\in\mathcal{O}_n$, it holds that
\[
\langle \lambda,\ \mu\rangle_{-} \le {\rm tr} (AXB^TX^T) \le \langle \lambda,\ \mu\rangle_{+}.
 \]
\end{thm}
The hidden convex reformulation of (O) reads as follows:
\begin{eqnarray}
&&\min_{X^TX =I_n} {\rm tr} (U{\rm diag}(\lambda)U^TXV{\rm diag}(\mu)V^TX^T) \nonumber \\
&=& \min_{X^TX =I_n} {\rm tr} ({\rm diag}(\lambda)U^TXV{\rm diag}(\mu)V^TX^TU)\nonumber\\
&=& \min_{Y^TY =I_n} {\rm tr} ({\rm diag}(\lambda)Y{\rm diag}(\mu)Y^T)\label{qap:1}\\
&\ge& \min_{\sum_{i=1}^ny_{ij}^2=\sum_{j=1}^ny_{ij}^2=1 ~\forall i,j}  \sum_{i=1}^n\sum_{j=1}^n\lambda_i\mu_jy_{ij}^2 \label{qap:2}\\
&=& \min_{\sum_{i=1}^nz_{ij}=\sum_{j=1}^nz_{ij}=1,z_{ij}\ge0 ~\forall i,j}  \sum_{i=1}^n\sum_{j=1}^n\lambda_i\mu_jz_{ij}, \label{qap:3}
\end{eqnarray}
where $Y=U^TXV$ is introduced in (\ref{qap:1}) so that $Y$ remains orthogonal, and in (\ref{qap:3}), $z_{ij}=y_{ij}^2$ $(i,j=1,\ldots,n)$ are introduced. Notice that the problem (\ref{qap:3}) is a linear assignment problem. Therefore, there is an optimal solution, denoted by $Z^*$,  lying at one of the vertices. That is, for all $i,j$, we have $z_{ij}^*\in \{0,1\}$. It follows that $y_{ij}^{*2}\in \{0,1\}$ and hence $y_{ij}^{*}\in \{0,1\}$, which further implies that  $Y^*\in {\rm\Pi}_n\subseteq \mathcal{O}_n$. Therefore, the inequality (\ref{qap:2}) holds as an equality.

This convexification  approach could be further extended to the trust-region type relaxation \cite{A99} by replacing (\ref{O1}) with
\[
\{X\in \Bbb R^{n\times n}:\ I_n\succeq X^TX\},
\]
and the enhanced version \cite{X11a} with (\ref{O1}) being replaced by
\[
\{X\in \Bbb R^{n\times n}:\ I_n\succeq X^TX,~{\rm tr}(X^TX)\ge m\},
\]
where $m$ is an integer between $0$ and $n$.

\section{Lagrangian dual and its variation}
In this section, we first show that strong Lagrangian duality could hold for a few nonconvex optimization problems or their special cases. Sometimes, in order to achieve strong duality,
the approach such as adding redundant constraints or making suitable transformation should be introduced in advance.

\subsection{Strong Lagrangian duality}

Let us begin with the generalized trust-region subproblem with interval bounds  \cite{P14}:
\begin{eqnarray*}
{\rm(QP)}~~&\min&  f_1(x) \\
&{\rm s.t.}& x\in\Omega:=\{x\in\Bbb R^n:\alpha\le f_2(x) \leq \beta\},
\end{eqnarray*}
where $f_i(x)=x^TA_ix+2a_i^Tx+c_i$ for $i=1,2$, $A_1,A_2\in \Bbb R^{n\times n}$ are symmetric matrices, $a_1,a_2\in \Bbb R^{n}$ and $c_1,c_2\in \Bbb R$. We make a further assumption:
\begin{ass}\label{as:1}
$A_2\neq 0$, $-\infty<\alpha\leq\beta<+\infty$ and there are $y,z\in \Bbb R^{n}$ such that either
 $\alpha<f_2(y)<\beta$ or $f_2(y)<\alpha=\beta<f_2(z)$ (i.e., the
primal Slater condition holds).
\end{ass}

A real application is the squared least squares model for the global positioning system (GPS) location \cite{B12}:
    \begin{equation}
    {\rm (SLS)}~~
    \min_{x\in\Bbb R^n,r\in\Bbb R} \left\{
    \sum_{i=1}^{m}\left(\|x-a_i\|^2-(r+d_i)^2\right)^2
    \right\},  \nonumber
    \end{equation}
  which  can be reformulated as a special case of (QP):
\[
    \min_{x\in\Bbb R^n,r,t\in\Bbb R}\left\{
 \sum_{i=1}^{m}\left(t-2a_i^Tx-2d_ir+a_i^Ta_i-d_i^2\right)^2:
 x^Tx-r^2=t
    \right\}.
    \]

A key technique to establish the strong duality is the following S-lemma with interval bounds \cite{W15}, which generalizes the classical S-lemma \cite{PT07} and the S-lemma with equality \cite{XSR16}.

\begin{thm}[S-lemma with interval bounds \cite{W15}] \label{Slemma1}
Under Assumption \ref{as:1}, the system $f_1(x)<0,~\alpha\leq f_2(x) \leq \beta$ is unsolvable if and only if there is a $\mu\in\Bbb R$ such that $f_1(x)+\mu_-(f_2(x)-\beta)+\mu_+(\alpha-f_2(x))\geq 0,~\forall x\in\Bbb R^n$, where $\mu_+=\max\{\mu,0\}$ and $\mu_-=\min\{\mu,0\}$.
\end{thm}
Then we show that the strong Lagrangian duality holds for (QP), which provides a hidden convexity of (QP) from the dual side.
\begin{eqnarray}
v{\rm (QP) }&=& \sup\left\{\lambda:\{x\in\Bbb R^n: f_1(x)-\lambda <0, \alpha\le f_2(x)\leq \beta\}=\emptyset\right\}\nonumber\\
&=&\sup\left\{\lambda:\exists \mu {\rm~such~that~} f_1(x)-\lambda  \right.\nonumber\\
&&~~~~~~~~\left.+\mu_-(f_2(x)-\beta)+\mu_+(\alpha-f_2(x))\geq 0, \forall x\in\Bbb R^n\right\}\label{e:0}\\
&=&\sup\left\{\lambda:
\left(\begin{array}{cc}A_1 -\mu A_2 & a_1 -\mu a_2\\
a_1^T -\mu a_2^T & c_1-\lambda  -\mu c_2-\mu_-\beta+\mu_+\alpha
\end{array}\right)\succeq 0\right\},\nonumber
\end{eqnarray}
where (\ref{e:0}) follows from Theorem \ref{Slemma1} and the last equality is based on the trivial observation
\begin{eqnarray*}
x^TAx+2a^Tx+c\geq 0, \forall x\in \Bbb R^n~\Longleftrightarrow~
\left(\begin{array}{cc}A & a\\
a^T & c\end{array}\right)\succeq 0.
\end{eqnarray*}
When $A_1$ and $A_2$ are diagonal or simultaneous
diagonalizable, (QP) admits a second-order conic reformulation \cite{BT14}, which seems to be easier to solve than the SDP.  Extensions to the general case can be found in \cite{J17}.

To study quadratic program with  two nonindependent  quadratic constraints, we need
the following general S-procedure due to Polyak \cite{Poly}.
\begin{thm} [Theorem 4.1, \cite{Poly}]
Let $n \ge3$ and $A_0,A_1, A_2$ be $n\times n$ real symmetric matrices. Suppose there are scalars $\mu_1,\mu_2$ and $\widetilde{x}\in\Bbb R^n$ such that
\begin{eqnarray}
&&\mu_1A_1+\mu_2A_2\succ 0,\label{S:1}\\
&&\widetilde{x}^TA_i\widetilde{x}<\alpha_i, ~i=1,2. \label{S:2}
\end{eqnarray}
Then, the system
\[
x^TA_0x<\alpha_0,~ x^TA_1x\le \alpha_1, ~  x^TA_2x\le \alpha_2
\]
has no solution if and only if there exist $\lambda_1\ge0$, $\lambda_2\ge0$ such that
\[
A_0+\lambda_1 A_1+\lambda_2A_2\succeq 0,~ \alpha_0+\lambda_1 \alpha_1+\lambda_2\alpha_2\le 0.
\]
\end{thm}
Consider the homogeneous nonconvex quadratic constrained quadratic programming:
\begin{eqnarray}
{\rm(QP2)}~~&\min&  x^TA_0x \label{qp2:1}\\
&{\rm s.t.}& x\in\Omega_2:=\{x\in\Bbb R^n:~x^TA_1x\le \alpha_1, ~  x^TA_2x\le \alpha_2\}.\label{qp2:2}
\end{eqnarray}
Under the assumption (\ref{S:1})-(\ref{S:2}), strong duality of (QP2) is verified as follows:
\begin{eqnarray}
v{\rm (QP2) }&=& \sup\left\{\lambda:\{x\in\Bbb R^n: x^TA_0x-\lambda <0, x^TA_1x\le \alpha_1, ~  x^TA_2x\le \alpha_2\}=\emptyset\right\}\nonumber\\
&=&\sup\left\{\lambda:
 A_0+\mu_1 A_1+\mu_2 A_2 \succeq 0,~  \lambda+\mu_1\alpha_1 +\mu_2\alpha_2 \le 0
\right\} \nonumber\\
&=&\sup\left\{-\mu_1\alpha_1 -\mu_2\alpha_2:
 A_0+\mu_1 A_1+\mu_2 A_2 \succeq 0
\right\}. \nonumber
\end{eqnarray}
Then we get a  semidefinite programming reformulation of (QP2). 

Strong duality of (QP2) can be extended to the problem with two equality constraints under the Slater's assumption for equalities.  As a small application, consider the binary quadratic programming problem
\begin{eqnarray*}
{\rm (BQP)}~~ \min \left\{ x^TQx:~ x\in \{-1,1\}^n\right\},
\end{eqnarray*}
which is NP-hard as it contains the classical Max-Cut  problem as a special case. The Lagrangian dual is given by
 \begin{eqnarray}
{\rm (D{\text-}BQP)}~~ \max \left\{ e^T\mu:~  Q-{\rm diag}(\mu) \succeq 0\right\}.
\end{eqnarray}
An optimal parametric Lagrangian dual approach is proposed in \cite{X15}, which 
implies  that
\[
v{\rm (BQP)}=v{\rm (D{\text-}BQP)}, ~n\le 2.
\]
Here, one can observe that  (BQP) with $n\le 2$ is just a special case of (QP2) with equality constraints and hence the strong duality holds true.

Note that (QP2) is homogeneous. However, for the inhomogeneous case, there could be a positive duality gap. A well-known example is the Celis-Dennis-Tapia (CDT) subproblem (see Section 5). But if the variables are in the complex field, strong duality holds again \cite{BE06}.

\subsection{Adding redundant constraints}

Consider the problem (O) (\ref{O10})-(\ref{O1}), the orthogonal relaxation of (QAP). Above we have shown that (O) enjoys hidden convexity. However, as pointed in
\cite{Z98}, the standard Lagrangian duality gap could be positive. Another counter example can be found in \cite{W02}.  Interestingly, the duality gap is closed by adding a  redundant constraint \cite{A00}. More precisely, (O) is equivalent to
\begin{eqnarray*}
{\rm (O2)}~~&\min& {\rm tr} (AXB^TX^T) \\
&{\rm s.t.}&   X^TX =I_n, \\
&&XX^T =I_n.
\end{eqnarray*}
The Lagrangian dual problem of (O2) is
\begin{eqnarray}
{\rm (DO2)}~~&\max& {\rm tr}(S+T) \nonumber\\
&{\rm s.t.}&  B\otimes A\succeq (I_{n}\otimes S)+(T
\otimes I_{n}),\nonumber\\
&&S=S^T\in\Bbb R^{n\times n},\ T=T^T\in\Bbb R^{n\times n}.\nonumber
\end{eqnarray}
It was shown in \cite{A00,W02} that  $v({\rm O2})=v({\rm DO2})$ and strong duality holds again.

Strong duality is achieved for the trust-region-type relaxation with additional constraints \cite{X11a}:
\begin{eqnarray}
&\min& {\rm tr} (AXB^TX^T) \nonumber \\
&{\rm s.t.}&    I_n\succeq X^TX, \nonumber \\
&& I_n\succeq XX^T, \label{trseig:1}\\
&&{\rm tr}(X^TX)\ge m,\nonumber
\end{eqnarray}
where $m\ge 0$ is an integer.   The special case $m=0$ was first studied in \cite{A99}.
Without the redundant constraint (\ref{trseig:1}), there could be a positive duality gap, see \cite{X11a}. For other applications of the similar approach, we refer to \cite{D11}.

The other example is the univariate polynomial optimization
\[
({\rm PO})~\min_{x\in\Bbb R}~P_{2k}(x):=a_{2k}x^{2k}+a_{2k-1}x^{2k-1}+\ldots+a_2x^2+a_1x+a_0,
\]
where $a_{2k}>0$ so that v(PO)$>-\infty$.
Let $x_i=x^i$. (PO) can be reformulated as the following quadratic program with many redundant quadratic constraints:
\begin{eqnarray}
 &\min& \sum_{i=0}^{k} a_ix_i+\sum_{i=k+1}^{2k} a_ix_ix_{i-k} \label{qcqp:1}\\
  &{\rm s.t.}& x_s-x_{s-t}x_t=0,~s=2,\ldots,k, ~1\le t\le  [(s+1)/2], \label{qcqp:2}\\
   & &x_kx_{s-k}-x_{s-t}x_{t}=0,~ k<s\le2k-2, ~1\le t\le  [(s+1)/2].\label{qcqp:3}
\end{eqnarray}
\begin{thm}[\cite{S98}]
  The Lagrangian dual value of the quadratic constrained quadratic program (\ref{qcqp:1})-(\ref{qcqp:3}) is equal to v$({\rm PO})$.
\end{thm}
The approach was further extended to the multidimensional  case with sum-of-square structure \cite{S98}.

\subsection{Scaled Lagrangian duality}
\subsubsection{Quadratic programming}
%

The standard quadratic programming is a nonconvex quadratic program over the standard simplex:
\begin{eqnarray}
{\rm (QPS)}~~ \min \left\{x^TQx:  ~x\in \Delta:=\{x\in  \Bbb R^n_+:~
e^Tx=1\}\right\},\nonumber
\end{eqnarray}
where $Q\in\Bbb R^{n\times n}$ and $Q\not\succeq 0$.  Any
general quadratic function $x^TQx+2b^Tx$ over $\Delta$ can be
homogenized by rewriting $b^Tx=x^T(eb^T)x$ and $x^Tb=x^T(be^T)x$ so
that $x^TQx+2b^Tx=x^T(Q+eb^T+be^T)x$. Problem (QPS) is NP-hard since  the maximum stability number can be reformulated in (QPS) \cite{M65}.

We can equivalently reformulate (QPS) as (see \cite{X13c})
\begin{eqnarray}
{\rm (QPS_1)}~~ \min \left\{ x^TQx:~   (e^Tx)^2=1,~x_ix_j\geq 0,~\forall i,j\right\}.\nonumber
\end{eqnarray}
The Lagrangian dual problem of ${\rm (QPS_1)}$ reads
\begin{eqnarray}
{\rm(D)}~~ \max \left\{ \lambda:~ Q-\lambda ee^T-S\succeq 0,~S\geq 0\right\}.\nonumber
\end{eqnarray}
The strong duality can be verified that \cite{A10}
\[
v({\rm QPS})=v({\rm D}),~ \forall ~n\le 4.
\]
It should be noted that there is a positive gap between (QPS) and its Lagrangian dual problem even when $n=5$.

\subsubsection{Fractional programming}

Reconsider the nonlinear convex fractional programming (NLF) (\ref{nlf:1})-(\ref{nlf:2}). As shown in  \cite{S76}, the standard Lagrangian dual is not tight for relaxing (NLF).  A natural way to define a dual problem with no duality gap is to write the Lagrangian dual of the equivalent convex reformulation (CLF) (\ref{nlf:3})-(\ref{nlf:5}), which reads as follows \cite{S74,S76}:
\begin{eqnarray}
\max_{\mu\ge 0} \left\{ \min_{g_i(x)\le 0,\forall i}~ \frac{g_0(x)+\sum_{i}\lambda_ig_i(x)}{h(x)}   \right\}.
\end{eqnarray}
In other words,    (NLF) with the following scaling
\begin{eqnarray*}
 \min \left\{  \frac{g_0(x)}{h(x)}:~
  \frac{g_i(x)}{h(x)}\le 0,~i=1,\ldots,m \right\}
\end{eqnarray*}
enjoys strong duality.

The assumption that $g_i(x)$ $(i=1,\ldots,m)$ are convex and $h(x)>0$ is concave  is not necessary in establishing strong duality of the above problem  (NLF).
Consider the quadratic fractional programming over an interval quadratic constraint:
\begin{eqnarray*}
{\rm(QPF)}~~\min \left\{ \frac{f_0(x)}{f_1(x)}:~x\in\Omega=\{x\in\Bbb R^n:\alpha\le f_2(x) \leq \beta\}\right\},
\end{eqnarray*}
where $f_i(x)$ ($i=0,1,2$) are quadratic functions.  Suppose Assumption \ref{as:1} holds.
To guarantee a well-defined (QPF), we further assume $\inf_{x\in\Omega} f_1(x)>0$. Special cases of the problem (QPF) were studied in \cite{B09,N16,X15}. In (QPF), $f_0(x)$ is nonconvex and $f_1(x)$ is nonconcave.

It was proved in \cite{Y17} that the equivalent reformulated optimization problem
\begin{eqnarray*}
 \min \left\{ \frac{f_0(x)}{f_1(x)}:~ \frac{\alpha}{f_1(x)}\le \frac{f_2(x)}{f_1(x)} \leq \frac{\beta}{f_1(x)} \right\}
\end{eqnarray*}
achieves zero Lagrangian duality gap. However, without this scaling, (QPF) could have a positive Lagrangian duality gap. Consider a special case of (QPF), the identical regularized total least squares  problem (TLS) \cite{B06}:
\[
{\rm(ITLS)}~~\min \left\{ \frac{\|Ax-b\|^2}{x^Tx+1}:~\alpha\le x^Tx \leq \beta\right\},
\]
where $A\in \Bbb R^{m\times n}$ ($m\ge n$) and  $0\le\alpha\leq\beta<+\infty$ is assumed.  The necessary and sufficient condition for the strong duality was established in \cite{Y17}.
\begin{thm}
The Lagrangian duality gap for $({\rm ITLS})$ is positive if and only if
\begin{eqnarray*}
&& \widetilde{A}:=A^TA-
\lambda_{\min}\left(\begin{array}{cc} A^TA & -A^Tb\\
-b^TA&b^Tb
\end{array}\right)\cdot I\succ  0,  \\
&&\alpha-b^TA  \widetilde{A}^{-2} Ab^T> 0.
\end{eqnarray*}
\end{thm}

\subsubsection{Orthogonal constrained linear optimization}
Replacing the objective function in (O) (\ref{O10})-(\ref{O1}) with a  linear one yields the following problem:
\begin{eqnarray}
({\rm LO})~~&\max& {\rm tr}(CX^T)\label{lo:1}\\
&{\rm s.t.}&X^TX=I_n.\label{lo:2}
\end{eqnarray}
It is easy to verify that the optimal value of (LO) is equal to the sum of all singular values of $C$.
However, there are examples such that the gap between (LO) and its Lagrangian dual is positive
\cite{W02}. To close the duality gap, (LO) is equivalently rewritten as
\begin{eqnarray*}
&\max& \frac{1}{2}{\rm tr}\left(\left(\begin{array}{cc}I_n&0\\0&0\end{array} \right)W\left(\begin{array}{cc}0&C^T\\C&0\end{array} \right)W^T \right)\\
&{\rm s.t.}&W^TW=I_n,~WW^T=I_n,
\end{eqnarray*}
where
$
W=\left(\begin{array}{cc}X&Y\\V&Z\end{array} \right).$ Then, we obtain (O2).

\section{Tight primal relaxation}
Primal relaxation could provide an alternative way to reveal the hidden convexity of the original nonconvex optimization problem.

\subsection{Totally unimodular}
Consider the integer programming problem
\[
\min_{Mx\le b,x\in Z^n}~ c^Tx,
\]
where $b\in Z^n$ and $Z^n$ is the field of integer vectors. The matrix $M$  is called totally unimodular if  each determinant  of any square submatrix is $0,1$ or $-1$. It follows that  the linear programming relaxation
\[
\min_{Mx\le b,x\in \Bbb R^n}~ c^Tx,
\]
has an integral optimal solution.  Instances of this class include  the min-cut problem, the linear assignment problem and so on.

\subsection{Orthogonal constrained problems}

We first consider a slightly generalized version of the linear problem (LO) (\ref{lo:1})-(\ref{lo:2}):
\begin{eqnarray*}
{\rm(GLO)}~~ \min \left\{ {\rm tr}(CX^T):~ X^TX=I_k, X\in\Bbb R^{n\times k}\right\},
\end{eqnarray*}
where $k\in \{1,\ldots,n\}$ and $C\in\Bbb R^{n\times k}$. (LO) (\ref{lo:1})-(\ref{lo:2}) exactly corresponds to the special case $k=n$.
(GLO) is trivially reformulated as the following convex programming:
\begin{eqnarray*}
 \min \left\{{\rm tr}(CX^T):~ X \in  {\rm conv}\left\{X\in\Bbb R^{n\times k}:~ X^TX=I_k\right\}\right\},
\end{eqnarray*}
which can be reduced to a semidefinite program as one can verify that
\begin{eqnarray*}
 {\rm conv}\left\{X\in\Bbb R^{n\times k}:~ X^TX=I_k\right\}&=&\left\{X\in\Bbb R^{n\times k}:~  I_k\succeq X^TX\right\}\\
 &=&\left\{X\in\Bbb R^{n\times k}:~ \left(\begin{array}{cc}I_n&X\\X^T&I_k\end{array} \right)\succeq 0  \right\}.
\end{eqnarray*}

Consider the relaxation problem (O) (\ref{O10})-(\ref{O1}). It is trivial to see that  (O) can be equivalently rewritten as a linear problem:
\begin{eqnarray*}
 \min \left\{ {\rm tr} (AY):~ Y\in \mathcal{O}(B):= {\rm conv}\{XBX^T:~ X\in
\mathcal{O}_n\}\right\}.
\end{eqnarray*}
The first computable representation of the convex set $\mathcal{O}(B)$  was established in \cite{X13} as follows:
\begin{thm}[\cite{X13}]\label{xia13}
\[
\mathcal{O}(B)= \left\{\sum_{i=1}^n\lambda_i(B)Y_i:
 \sum_{i=1}^nY_i=I_n, {\rm tr} (Y_i)=1, Y_i\succeq 0, Y_i=Y_i^T \in
\Bbb R^{n\times n}, \forall i\right\}.
\]
\end{thm}
Then, we obtain a semifinite programming reformulation of (O).

\subsection{Trust-region subproblem and extensions}
The trust-region subproblem (TRS)  (\ref{trs}) admits the following primal convex quadratic optimization reformulation \cite{F96}
\begin{eqnarray}
&\min& x^T(A-\lambda_{\min}(A))x+b^Tx +\delta\lambda_{\min}(A) \label{ct:1} \\
&{\rm s.t.}& x^Tx\le \delta.  \label{ct:2}
\end{eqnarray}
When extended to the two-sided trust-region subproblem
\begin{eqnarray*}
\min \left\{ x^TA x+b^Tx:~ \alpha\le x^Tx\le \delta\right\},
\end{eqnarray*}
the following equivalent convex reformulation was established in \cite{W17}:
\begin{eqnarray*}
&\min& x^T(A-\lambda_{\min}(A))x+b^Tx + \lambda_{\min}(A)\cdot t \\
&{\rm s.t.}& x^Tx\le t,~ \alpha\le t\le \delta.
\end{eqnarray*}
Applying Nesterov¡¯s accelerated gradient descent algorithm for solving (\ref{ct:1})-(\ref{ct:2}) yields  a  linear-time algorithm \cite{N17,W17}, whose worst-case complexity is less than the previously existing algorithm \cite{H16}.
Recently, the extended (TRS) where (\ref{ct:2}) is replaced by a general quadratic constraint  is convexified in  \cite{Y18}.

Besides (TRS), some other optimization problems over the unit-ball  also have hidden convexity.  For example,  the ball-constrained  weighted maximin dispersion problem \cite{Ha13}
\begin{eqnarray*}
{\rm (DP)}~~\max_{\|x\|\le 1}\left\{\min_{i=1,\ldots,m} \omega_i\|x-x^i\|^2\right\},
\end{eqnarray*}
where $x^1,\ldots,x^m\in \Bbb R^n$ are given $m$ points and $\omega_i>0$ for $i=1,\ldots,m$.  Problem (DP) is NP-hard \cite{W16}. However, when $m\le n$, (DP) is equivalent to
the following  second-order cone programming  problem \cite{W16}:
\begin{eqnarray*}
&\max&~ \zeta\\
&{\rm s.t.}& \omega_i(1-2(x^i)^Tx+\|x^i\|^2)\geq \zeta,~i=1,\ldots,m,\\
&& \|x\|\leq 1.
\end{eqnarray*}
The other example comes from the identical Tikhonov regularized total least squares \cite{BB06,B06}
\begin{equation}
{\rm (TI)}~~\min_{x\in \Bbb R^n} \left\{\frac{\|Ax-b\|^2}{\|x\|^2+1}+\rho\|x\|^2 = \frac{\|Ax-b\|^2+\rho\|x\|^4+\rho\|x\|^2}{\|x\|^2+1} \triangleq  \frac{f(x)}{g(x)}\right\}.
\label{TI}
\end{equation}
According to the S-lemma with equality \cite{XSR16}, a special case of Theorem \ref{Slemma1} with the setting $\alpha=\beta$,  we have
\begin{eqnarray}
v{\rm (TI)}=&&\max\{t: \min_{x\in \Bbb R^{n}} \{f(x)-tg(x)\}\ge 0\}\nonumber\\
=&&\max\{t: \{x\in \Bbb R^n:~f(x)-tg(x)<0\}=\emptyset \}\nonumber\\
=&&\max\{t: \{(x;s): \|Ax-b\|^2+\rho s^2+\rho s-t (s+1)<0,\|x\|^2=s\}=\emptyset \}\nonumber\\
=&& \max \{t:\|Ax-b\|^2+\rho s^2+\rho s-t (s+1)+\mu(\|x\|^2-s)\geq 0,\forall (x;s) \}\nonumber\\
=&&\max\left\{t: \left(\begin{array}{ccc}A^TA+\mu I &0&-A^Tb\\
0&\rho &\frac{\rho-t-\mu}{2}\\
-b^TA&\frac{\rho-t-\mu}{2}&b^Tb-t
\end{array}\right)
\succeq 0\right\}. \nonumber
\end{eqnarray}
The above SDP reformulation can be equivalently reduced to finding the unique zero point of a smooth, strictly decreasing and convex univariate function in terms of $\mu$. For more details, see \cite{Y17b}.

Moreover,
(TRS) (\ref{trs}) itself has the following equivalent semidefinite programming reformulation \cite{R97}:
\begin{eqnarray}
&\min &  {\rm tr}(AX)+c^Tx \nonumber\\
&{\rm s.t.}&{\rm tr}(X) \le \delta, ~ X\succeq xx^T.\nonumber
\end{eqnarray}
The basic idea behind this equivalence is Pataki's theorem \cite{P98}.
\begin{thm}[\cite{P98}]\label{Pa}
Suppose the SDP problem
\begin{eqnarray*}
{\rm(SDP)}&\min& {\rm tr}(A_0X)\\
&&{\rm tr}(A_iX)\le b_i,i\in I,\\
&&{\rm tr}(A_iX)=b_i,i\in E,\\
&&X\succeq 0
\end{eqnarray*}
has a solution and assume
\[
|I| + |E| \le (r+2)(r+1)/2-1.
\]
Then, for (SDP), there is an optimal solution $X^*$ such that rank($X^*$)$\le r$.
\end{thm}
Theorem \ref{Pa} implies that (QP2) (\ref{qp2:1})-(\ref{qp2:2}) admits a tight primal SDP relaxation.

(TRS) was further generalized by adding several linear cuts:
\begin{eqnarray}
({\rm T_m})~~  &\min & x^TAx+c^Tx \nonumber\\
&{\rm s.t.}&x^Tx\le \delta,\label{Tm2}\\
&&a_i^Tx\le b_i,~i=1,\ldots,m.\label{Tm3}
\end{eqnarray}
When $m=1$, the hidden convex reformulation of $({\rm T_1})$ was given in \cite{S03,Y03}:
\begin{eqnarray}
&\min &  {\rm tr}(AX)+c^Tx \nonumber\\
&{\rm s.t.}& \|b_1x-Xa_1\|\le \sqrt{\delta}(b_1-a_1^Tx),\label{tr1}\\
&&{\rm tr}(AX) \le \delta, ~ X\succeq xx^T,\nonumber
\end{eqnarray}
where the additional second-order constraint (\ref{tr1}) was obtained by linearizing the valid constraint:
\[
\|(b_1-a_1^Tx)x\|=(b_1-a_1^Tx)\|x\|\le \sqrt{\delta}(b_1-a_1^Tx).
\]
Further extension such as letting $\delta$ be an additional variable is studied in \cite{G14}.

Suppose $m=2$ and the two linear cuts are parallel, without loss of generality, we assume the constraints (\ref{Tm3}) are
\[
l \le a_1^Tx \le u.
\]
Then, $({\rm T_2})$ is equivalent to the following SOC-SDP problem \cite{B13}:
\begin{eqnarray}
&\min &  {\rm tr}(AX)+c^Tx \nonumber\\
&{\rm s.t.}& a_1^TXa_1+lu\le (l+u)a_1^Tx, \label{tr2}\\
&&\|ux-Xa_1\|\le \sqrt{\delta}(u-a_1^Tx),\nonumber \\
&&\|lx-Xa_1\|\le \sqrt{\delta}(a_1^Tx-l),\nonumber \\
&&{\rm tr}(AX) \le \delta, ~ X\succeq xx^T,\nonumber
\end{eqnarray}
where (\ref{tr2}) corresponds to linearizing the valid constraint:
\[
(a_1^Tx- u)(a_1^Tx- l)=a_1^Txx^Ta_1-(l+u)a_1^Tx+lu \le 0.
\]
Generally, when all the linear cuts (\ref{Tm3}) are non-intersecting in the ball (\ref{Tm2}), $({\rm T_m})$ has the following SOC-SDP reformulation \cite{B15}:
\begin{eqnarray}
&\min &  {\rm tr}(AX)+c^Tx \nonumber\\
&{\rm s.t.}& b_ib_j-b_ja_i^Tx-b_ia_j^Tx +a_i^TXa_j\ge 0,~i,j=1,\ldots,m,i<j, \nonumber\\
&&\|b_ix-Xa_i\|\le \sqrt{\delta}(b_i-a_i^Tx),~i=1,\ldots,m,\nonumber \\
&&{\rm tr}(AX) \le \delta, ~ X\succeq xx^T.\nonumber
\end{eqnarray}
Besides, for the intersection of an ellipsoid and a split disjunction, the convex hull is recently shown to be second-order-cone representable \cite{B17}.

\subsection{Quadratic matrix programming}
Consider a class of quadratic matrix programming \cite{Be06}:
\begin{eqnarray}
&\min_{X\in\Bbb R^{n\times r}}&  {\rm tr}(X^TA_0X)+2{\rm tr}(V^TB_0^TX)+c_i \nonumber\\
&{\rm s.t.}&{\rm tr}(X^TA_iX)+2{\rm tr}(V^TB_i^TX)+c_i\le \alpha_i,~i\in I, \nonumber\\
&&{\rm tr}(X^TA_iX)+2{\rm tr}(V^TB_i^TX)+c_i= \alpha_i,~i\in E, \nonumber
\end{eqnarray}
where $A_i=A_i^T\in\Bbb R^{n\times n}$, $B_i\in\Bbb R^{n\times s}$ ($i\in\{0\}\cup I\cup E$) and $V\in R^{s\times r}$ with $s\le r$. Applications of this model include robust least squares and the sphere-packing problem.

For $i\in\{0\}\cup I\cup E$, define
\[
M_i=\left(\begin{array}{cc}A_i&B_i\\B_i^T&\frac{c_i}{{\rm tr}(VV^T)I_s}\end{array} \right).
\]
Relaxing $\left(\begin{array}{c}X\\V\end{array} \right)(X^T~~V^T)$ to $Z$ yields the following SDP:
\begin{eqnarray}
&\min_{X\in\Bbb R^{n\times r}}&  {\rm tr}(M_0Z) \nonumber\\
&{\rm s.t.}&{\rm tr}(M_iZ)\le \alpha_i,~i\in I, \nonumber\\
&&{\rm tr}(M_iZ)= \alpha_i,~i\in E, \nonumber\\
&&Z\succeq 0,\nonumber\\
&&Z_{n+i,n+j}=(VV^T)_{i,j},~i,j=1,\ldots,s.\nonumber
\end{eqnarray}
As an application of Theorem \ref{Pa}, it was shown in \cite{BD12} that
the above primal  SDP relaxation is tight if its optimal value is attainable and either $n+s\le r$ or
$|I|+|E|\le \frac{(r+2)(r+1)}{2}-\frac{s(s+1)}{2}-1$.

For more nonconvex instances with tight SDP relaxation, we refer to \cite{L16,L13,L97,Y18} and references therein.

\section{Open problems}

We conjecture that any polynomially solved optimization problem has an equivalent hidden convex reformulation. Nevertheless, we are expected to
find the hidden convex reformulations of the following ten special nonconvex optimization problems in the near future.

\noindent{\bf Open Problem 1.} The  Celis-Dennis-Tapia (CDT) subproblem
\begin{eqnarray*}
({\rm CDT})&\min&  x^TQ_0x + 2q_0^Tx, \\
 &{\rm s.t.} &  x^TQ_ix + 2q_i^Tx + {\gamma _i} \le 0,~i = 1,2,
\end{eqnarray*}
where $Q_1\succ 0$. Recently, the polynomial-time solvability of (CDT) has been proved in \cite{B16,CL2016,S16}. Hidden convex reformulation even for the diagonal (CDT) (i.e., $Q_0,Q_1$ and $Q_2$ are all diagonal matrices) remains unknown.

\noindent{\bf Open Problem 2.} Finding the hidden convex reformulation of the extended trust-region subproblem with \textit{any} fixed number of linear cuts
\begin{eqnarray*}
&\min&  x^TQ_0x + 2q_0^Tx + {\gamma _0}, \\
 &{\rm s.t.} &  x^Tx\le 1,\\
 && a_i^Tx\le b_i,~i=1,\ldots,k,
\end{eqnarray*}
which is polynomially solved when $k$ is fixed \cite{B14,H13}.

\noindent{\bf Open Problem 3.}  The regularized version of the trust-region subproblem with $k$ linear cuts reads as follows:
\begin{eqnarray*}
&\min&  x^TQ_0x + 2q_0^Tx + \|x\|^p, \\
 &{\rm s.t.} &  a_i^Tx\le b_i,~i=1,\ldots,k,
\end{eqnarray*}
where $p>2$ is a parameter. When $k=0$, there is at most one local non-global minimizer \cite{H17}. Hidden convexity with any fixed $k\ge 1$ is unknown.

\noindent{\bf Open Problem 4.} The unbalanced orthogonal Procrustes problem
\begin{eqnarray*}
\min_{X^TX=I_k} ~\|AX-B\|_F^2.
\end{eqnarray*}
Notice that when $k=1$, it reduces to trust-region subproblem. And when $k=n$, it reduces to the balanced case with an explicit solution \cite{S66} and admits a quadratic matrix programming \cite{Be06}. Hidden convex reformulation for the special cases $k=2$ and $k=n-1$ are expected.

\noindent{\bf Open Problem 5.} Finding the hidden convex reformulation of the Tikhonov regularized  total least squares problem \cite{BB06}
\begin{eqnarray*}
\min  \frac{\|Ax-b\|^2}{x^Tx+1}+\rho \|Lx\|^2,
\end{eqnarray*}
which is polynomially solved as it can be reformulated as a special case of the generalized (CDT). The special case $L=I$ is settled in \cite{Y17b}.

\noindent{\bf Open Problem 6.} Finding the hidden convex reformulation of the sum of a generalized Rayleigh quotient and a quadratic form on the unit sphere
\begin{eqnarray*}
\min_{\|x\|=1}  \frac{x^TAx}{x^TBx}+x^TCx,
\end{eqnarray*}
where $B\succ0$. This problem is a generalization of the Rayleigh quotient optimization problem. It was raised in \cite{Z14} with applications in the downlink of a multi-user MIMO system and  the sparse Fisher discriminant analysis in pattern recognition.

\noindent{\bf Open Problem 7.}  Let $A,B$ be $n\times n$ positive definite matrices and $0\neq b\in\Bbb R^n$. The special unconstrained quartic minimization
\begin{eqnarray*}
\min_{x\in\Bbb R^n}   (x^TAx)(x^TBx)+b^Tx,
\end{eqnarray*}
has special application in the
Legendre-Fenchel conjugate of the product of two positive definite quadratic functions \cite{H07}. Assuming the objective function being convex, it has been solved in \cite{Z10}. It is further shown in \cite{X13a} that the convexity assumption could be removed.  However, hidden convexity of this problem remains unknown.

\noindent{\bf Open Problem 8.}  Let $A_k\in\Bbb R^{n\times n}$ be symmetric for $k=1,2,\ldots,m$. The ball-constrained quartic minimization
\begin{eqnarray*}
\min_{\|x\|=1}  ~\sum_{k=1}^m (x^TA_kx)^2
\end{eqnarray*}
is generally NP-hard \cite{N03}.  Hidden convexity even for the special case $k=2$ remains unknown.

\noindent{\bf Open Problem 9.}  Let $A_k\in\Bbb R^{n\times n}$ be symmetric for $k=1,2,\ldots,m$. The optimal value of
\[
\max_{\|x\|=1} \lambda_{\max}\left([A_1x~\ldots~A_mx]^T[A_1x~\ldots~A_mx]\right).
\]
plays a great role in the local convex analysis for quadratic
transformations \cite{X14}. When $m=1$, the optimization problem reduces to the maximal eigenvalue of $A_1^TA_1$. Hidden convexity for fixed $m$ (say $m=2$) is unknown.

\noindent{\bf Open Problem 10.} Let
$\Omega=\{x\in \Bbb R^n: \|x-a_i\|^2\leq r_i^2,~i=1,\ldots, p\}$ be the intersection of $p$ balls.
Finding the Chebyshev center of $\Omega$ is modeled as
\begin{equation}
{\rm (CC_B)}~~\min_{z}\max_{x\in \Omega}\|x-z\|^2.
\end{equation}
Geometrically, $({\rm CC_B})$ is to find the smallest ball enclosing $\Omega$. When $p\le n$, $({\rm CC_B})$ admits a standard quadratic programming representation \cite{Be07,Be09}.  Hidden convexity for fixed $n$ or $p=n+1$ remains unknown.

\end{document}